\documentclass[12pt]{iopart}

\expandafter\let\csname equation*\endcsname\relax 
\expandafter\let\csname endequation*\endcsname\relax 

\usepackage{amsmath}
\usepackage{amsfonts}
\usepackage{graphicx}
\usepackage{subfigure}
\usepackage{soul}
\usepackage{color}

\begin{document}

	\title{The role of detachment of in-links in scale-free networks}
		
	\author{P Lansky$^1$, F Polito$^2$, L Sacerdote$^2$}

	\address{$^1$ Institute of Physiology, Academy of Sciences of the Czech Republic,
			Videnska 1083, Praha 4, Czech Republic}
	\address{$^2$ Dipartimento di Matematica, Universit\`{a} degli Studi di Torino,
			Via Carlo Alberto 10, 10123 Torino, Italy}

	\begin{abstract}
		Real-world networks may exhibit detachment phenomenon determined by the cancelling of previously existing connections. We
		discuss a tractable extension of Yule model to account for this feature.
		Analytical results are derived and discussed 
		both asymptotically and for a finite number of links.
		Comparison with the original model is performed in the supercritical case.
		The first-order asymptotic tail behavior of the two models is similar but differences arise in the second-order term.
		We explicitly refer to World Wide Web modeling and we show the agreement of the proposed model on very recent data. However,
		other possible network applications are also mentioned.		
	\end{abstract}
	

	\section{Introduction}
	
		Large scale hypertexts, such as the World Wide Web, are well-known examples of scale-free networks in the real-world \cite{pietro}
		and the Barab\'asi--Albert model \cite{ba} is certainly one of their most investigated descriptors.
		When the interest focuses on the popularity of a webpage, its connectivity (i.e.\ the number of its in-links)
		is usually considered. On the other hand, the knowledge of in- and out-links becomes necessary to deal with the topology of the network.
		Both approaches exist in the literature for the so called  \textit{preferential attachment model} \cite{ba}, in which
		new webpages are constantly added to the network
		with a given probability and the number of in-links of each webpage grows proportionally to the number of in-links already possessed
		by that webpage. It turns out that with these assumptions, the limiting distribution of the number
		of in-links for a webpage chosen at random exhibits a power-law tail (see \cite{ba} for a more in-depth explanation).
		It is pointed out in \cite{born} that the Barab\'asi--Albert model 
		for in-links
		can be classified as a specific case of the so-called
		Yule--Simon model \cite{simon} which in turn is related to the well-known Yule model for macroevolution \cite{yule}
		(do not confuse it with the \emph{Yule process} which is a homogeneous linear birth process).
		Specifically, Yule--Simon and Yule models are asymptotically equivalent \cite{willis}.
		The Yule model appears recurrently with different names in relation to different applications \cite{willis}.
		Despite the long history of the Yule model, it maintains a relevant role in a variety of applications
		\cite{doro}; variants of it have been studied in the past \cite{doro2,new} and are part of very recent research \cite{pnas,plos}.
				
		The aim of this paper is to study formally the consequences of detachment of in-links.
		This is accomplished by considering the Yule model and replacing the linear birth process governing
		the growth of in-links with a linear birth-death process (following \cite{reed}, where the model was first introduced
		within the context of macroevolution).
		It turns out that the introduction of the possibility of detachment of in-links in the Yule model still leads to an
		analytically tractable model and thus still permits to obtain exact results.
	
		In Section \ref{inaina} we determine the analytical form of the distribution of the number of in-links to a webpage.
		Our results include the probabilities
		of any finite number of in-links not being limited to the asymptotic
		limits, as happens in the literature.
		In fact, the analysis of the role of the detachment becomes possible for any range of the number of in-links.
		The availability of the in-links distribution allows us the comparison of its features with those of the World Wide Web in 2012
		(Section \ref{ddataa}).
		Further, the moments are calculated explicitly and their expressions are
		presented in Section \ref{mome}.
		In Section \ref{fuego} we compare the generalized model and the classical Yule model.

	\section{Classical model}

		The construction based on Yule model \cite{yule} proceeds as follows.
		Consider an initial webpage present at time $t=0$ with a single in-link.
		The number of in-links increases proportionally to the number of existing in-links with
		an in-link rate constant $\lambda>0$.			
		The number of webpages develops independently as a linear birth process
		of webpage rate constant $\beta>0$.
						
		\begin{figure}
			\centering
			\includegraphics[scale=.6]{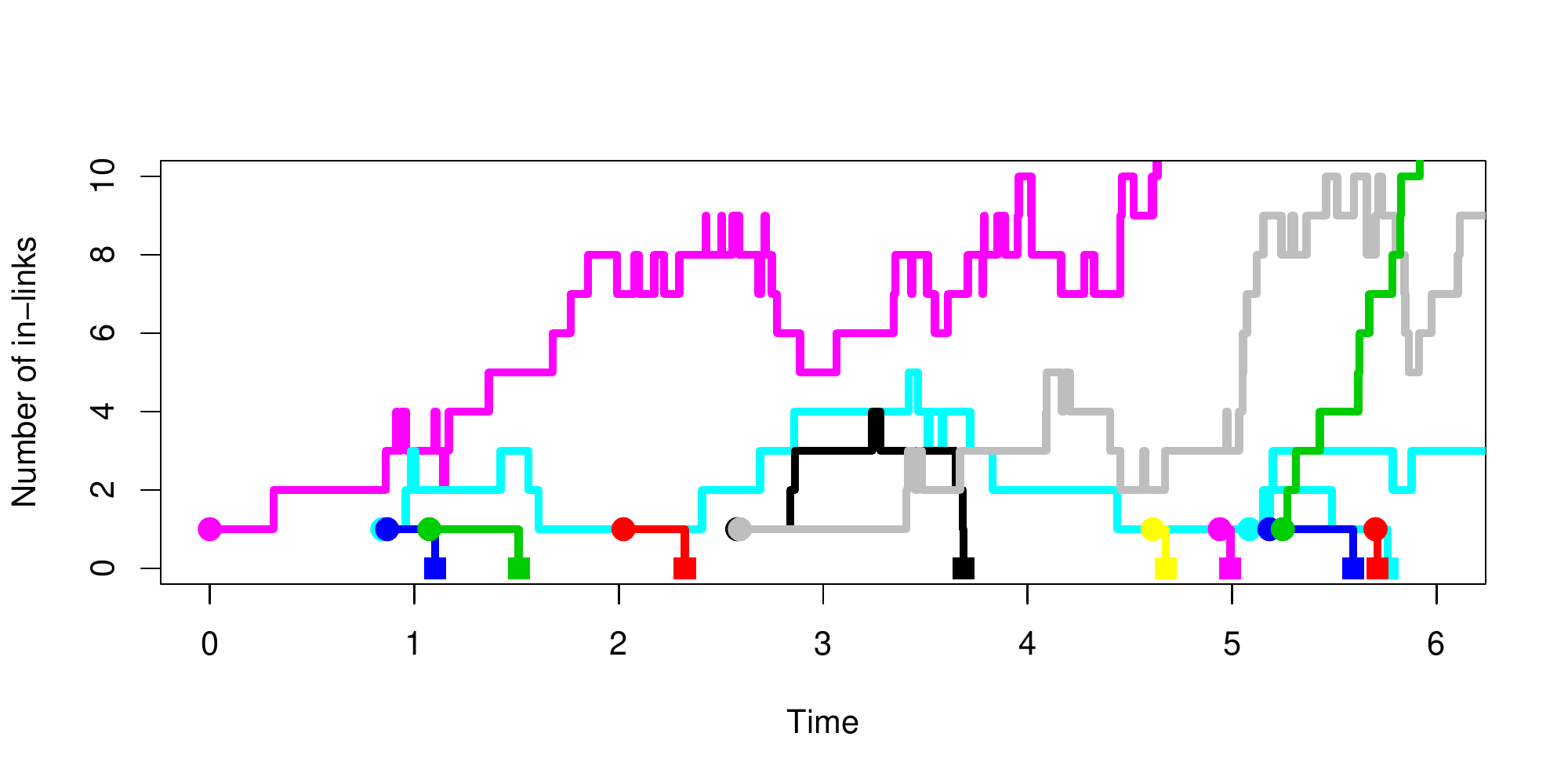}
			\caption{\label{figure1}An example of a realization of the generalized Yule model with
				$\beta=0.1$, $\lambda=1.1$, $\mu=1$ (supercritical case). The evolution of the number of in-links in different webpages is
				highlighted with different colors. The instants when new webpages are introduced are indicated with colored dots
				while the moments when webpages disappear due to the removal of the last in-link
				are denoted by colored squares.}
		\end{figure}
		Let $N_\beta(t)$ and $N_\lambda(t)$ be the counting processes reflecting the number of webpages and of in-links of a given webpage,
		respectively, observed at time $t$.
		The assumption that the process $N_\beta(t)$ starts
		at time zero with one initial webpage simplifies the notation but
		has no effects on the asymptotic properties of the results.
		Recall that the probability distribution $\mathbb{P} ( N_\lambda(t) = n )$, $n \ge 1$, of
		a linear birth process is geometric \cite{bailey},
		\begin{align}
			\label{bai}
			\mathbb{P} ( N_\lambda(t) = n ) = e^{-\lambda t} (1-e^{-\lambda t})^{n-1}, \quad t \ge 0.
		\end{align}
		In the above formula, as removal of in-links is not permitted, $n$ cannot reach zero.
		
		Referring to
		the results obtained in \cite{neuts,crump,feigin,puri}, by conditioning
		on the number of webpages present at time $t$, we obtain that the
		random instants at which appearance of new webpages occurs are distributed as the order statistics
		of independent and identically distributed random variables $T$ with probability distribution function
		\begin{align}
			\label{distr}
			\mathbb{P} ( T \le \tau ) = \frac{e^{\beta \tau}-1}{e^{\beta t}-1},  \qquad 0 \le \tau \le t.
		\end{align}
		Let $\mathcal{N}^Y_t$ denote the size of a webpage chosen uniformly at random at time $t$ and call $\mathcal{N}^Y
		= \lim_{t \to \infty}\mathcal{N}^Y_t$.
		Plainly,
		\begin{align*}
			\mathbb{P} (\mathcal{N}^Y_t = n ) = \mathbb{E}_T \mathbb{P} ( N_\lambda(t) = n | N_\lambda(T) = 1 ),
		\end{align*}
		meaning that the selected webpage was introduced at a random time $T$ with
		distribution \eqref{distr}. Note that the age $t-T$ of the randomly selected webpage appeared at time $T$ is
		distributed as a truncated exponential random variable.
		Hence, using \eqref{bai} we have
		\begin{align}
			\label{lenovo}
			\mathbb{P} (\mathcal{N}^Y_t = n )
			= \frac{\beta}{1-e^{-\beta t}} \int_0^t e^{-\beta y} e^{-\lambda y} (1-e^{-\lambda y})^{n-1} \mathrm dy,
		\end{align}
		and in the limit for $t \to \infty$ we obtain
		\begin{align}
			\label{two}
			\mathbb{P} ( \mathcal{N}^Y = n ) & = \int_0^\infty \beta e^{-\beta y} e^{-\lambda y}
			(1-e^{-\lambda y})^{n-1} \mathrm dy.
		\end{align}
		Expression \eqref{two} is simple
		to evaluate (see for example \cite{yule}),
		\begin{align}
			\label{yule}
			\mathbb{P}(\mathcal{N}^Y=n) = \frac{\beta}{\lambda} \frac{\Gamma(n)\Gamma
			\left( 1+ \frac{\beta}{\lambda} \right)}{\Gamma \left( n+1+\frac{\beta}{\lambda} \right)}.
		\end{align}
		Distribution \eqref{yule} is known as Yule--Simon distribution \cite{bala}.

	\section{Generalized model: basic properties}
		
		\label{inaina}
		As stressed above, we are concerned with the limiting behavior of a stochastic process in which the
		appearance (observation) of new webpages (each with only one in-link) is controlled by a linear birth process $N_\beta(t)$
		of intensity $\beta>0$. 
		Now the evolution of the number of in-links for each webpage is governed by
		independent (of one another and of $N_\beta(t)$) linear birth-death processes $N_{\lambda,\mu}(t)$ with birth and death intensities
		$\lambda>0$ and $\mu>0$, respectively (see also the description of this model in
		\cite{reed}). A schematic example of the model can be seen in Fig.\ \ref{figure1}.
		Note how, with that particular choice of the parameters, many webpages may lose all of their links in a short time.
		However, a few of them develop fast gaining a considerable number of in-links in a relatively short timespan.
											
		To derive the asymptotic properties of the model,
		we proceed by following the same lines leading to formulae \eqref{two} and \eqref{yule}. To compute the distribution
		of the number of in-links $\mathcal{N}$ of a webpage chosen uniformly at random,
		we first write the analogue of \eqref{lenovo} randomizing by the time $T$ of first
		appearance of a webpage and then we let $t \to \infty$.
		Since in this generalized model	$\mathcal{N}$ can be zero with a positive probability,
		in the following, we distinguish the case $\mathbb{P}(\mathcal{N}=0)$ from $\mathbb{P}(\mathcal{N}=n)$, $n \ge 1$.
		For the sake of brevity we will omit $n \ge 1$ in formulae.

		As we are dealing with birth-death processes for the evolution of the number of in-links,
		we must separate the three cases $\lambda=\mu$ (critical), $\lambda<\mu$ (subcritical), and $\lambda>\mu$ (supercritical)
		in which the processes behave rather differently (e.g.\ \cite{bailey}, Section 8.6).

		\subsection{Critical regime}

\begin{figure}
				\centering
				\subfigure[]{\includegraphics[scale=0.75]{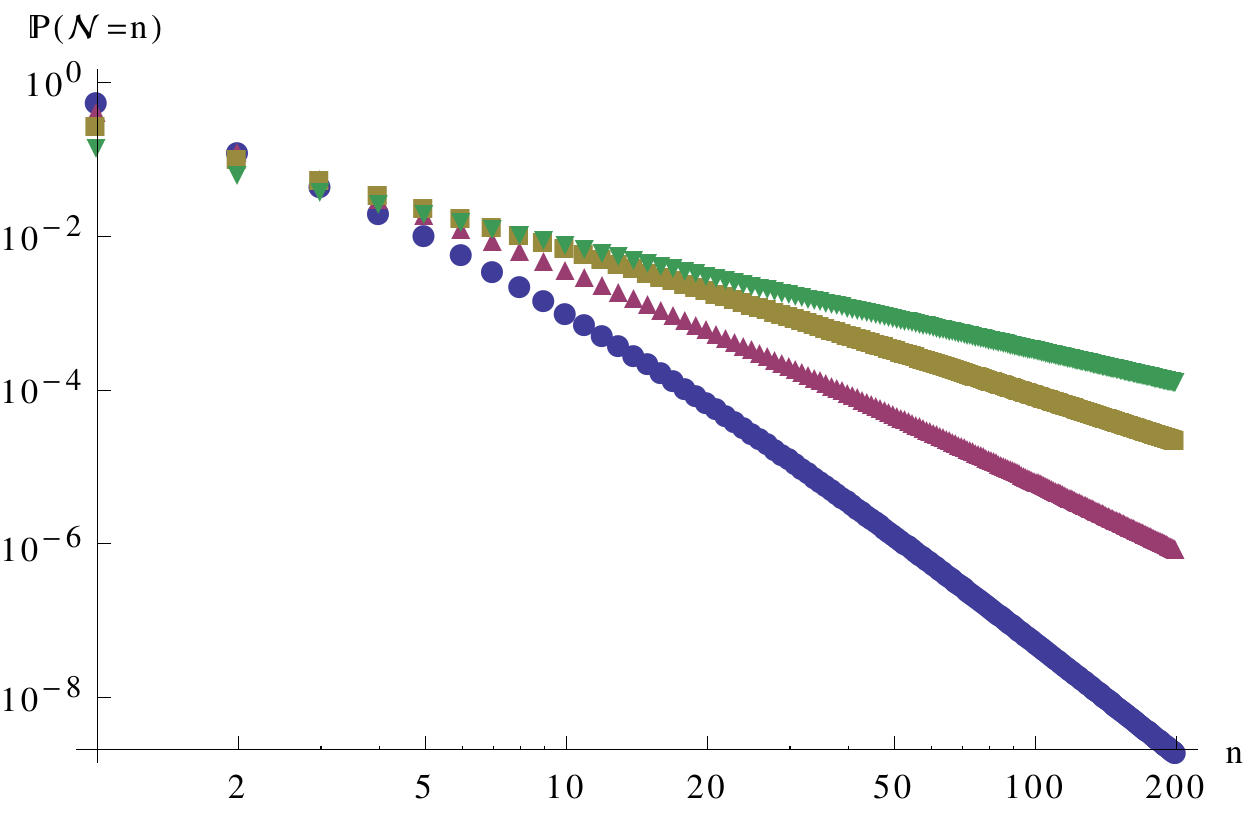}}
				\subfigure[]{\includegraphics[scale=0.75]{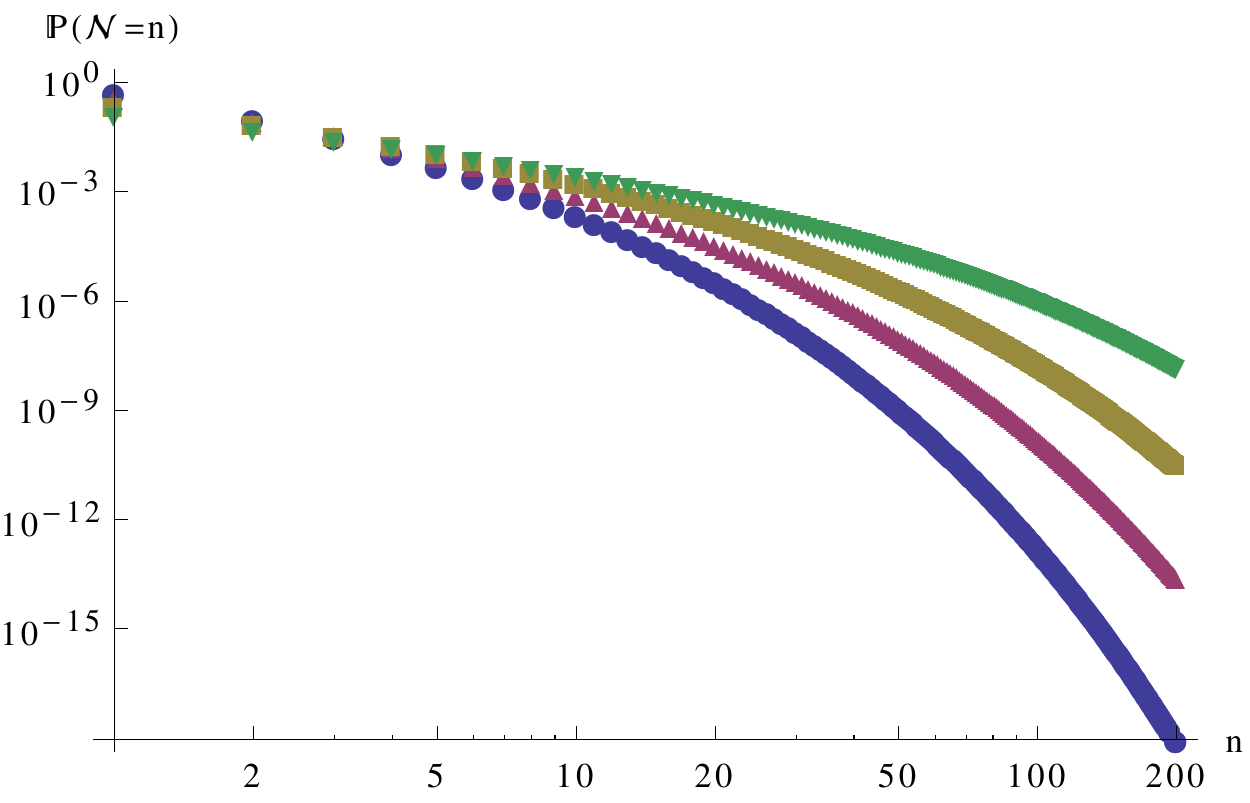}}
				\subfigure[]{\includegraphics[scale=0.75]{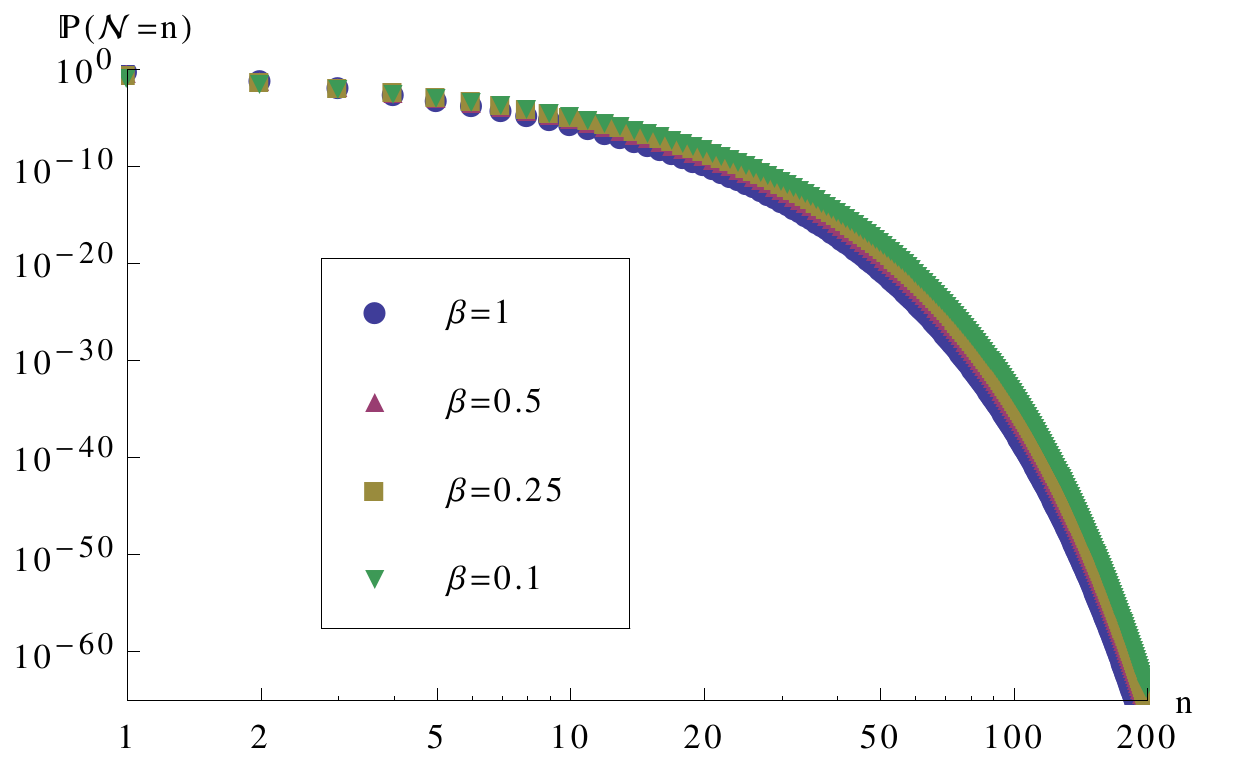}}
				\caption{\label{figure2b}
				Distribution of the number of in-links for different values of the webpage rate constant, $\beta=\{1,0.25,0.5,0.1\}$.
				(a) Supercritical case, $\lambda=1/2$, $\mu=1/4$; the tails decay as power-laws (see the discussion in Section \ref{fuego}).
				(b) Critical case, $\lambda=\mu=1/2$;
				we can see that the tails decay faster than a power-law.
				(c) Subcritical case, $\lambda=1/4$, $\mu=1/2$; note that
				the webpage rate constant $\beta$ plays almost no role.}
			\end{figure}

			Let us start with the case $\lambda=\mu$. Considering that for $n \ge 0$,
			\begin{align}
				\label{su}
				\mathbb{P}(\mathcal{N}=n) = \int_0^\infty \beta e^{-\beta t} \mathbb{P} (N_{\lambda,\mu}(t) = n) \, \mathrm dt,
			\end{align}
			we have
			\begin{align*}
				\mathbb{P}(\mathcal{N}=0) = \int_0^\infty \beta e^{-\beta t} \frac{\lambda t}{1+\lambda t} \mathrm dt
			\end{align*}
			and
			\begin{align*}
				\mathbb{P}(\mathcal{N}=n)
				= \int_0^\infty \beta e^{-\beta t} \frac{(\lambda t)^{n-1}}{(1+\lambda t)^{n+1}} \, \mathrm dt.
			\end{align*}
			Exploiting the integral representation of the confluent hypergeometric function (see \cite{abramowitz}, Section 13)
			\begin{align}
				\label{pdf}
				U(a,b,z) = \frac{1}{\Gamma(a)} \int_0^\infty e^{-zy} y^{a-1} (1+y)^{b-a-1} \, \mathrm dy
			\end{align}
			and the transformation formula 13.1.29 of \cite{abramowitz},
			the above expressions can be written as
			\begin{align}
				\mathbb{P}(\mathcal{N}=0) = U(1,0,\beta/\lambda),
			\end{align}
			and
			\begin{align}
				\label{ext}
				\mathbb{P}(\mathcal{N}=n) =
				(\beta/\lambda) \, \Gamma(n) \, U(n,0,\beta/\lambda).
			\end{align}
			Fig.\ \ref{figure2b}(b) shows the probability mass function of $\mathcal{N}$ in the specific critical case $\lambda=\mu=1/2$.
			Notice the fast decay of the tail.
	
		\subsection{Non-critical regimes}
				
\begin{figure}
			\centering
			\subfigure[]{\includegraphics[scale=0.75]{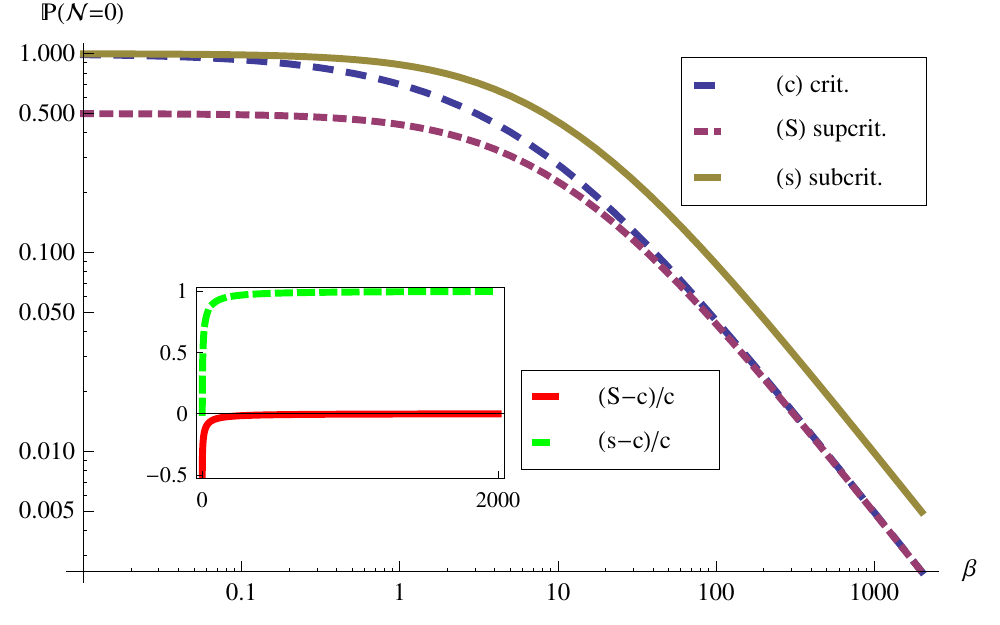}}
			\subfigure[]{\includegraphics[scale=0.85]{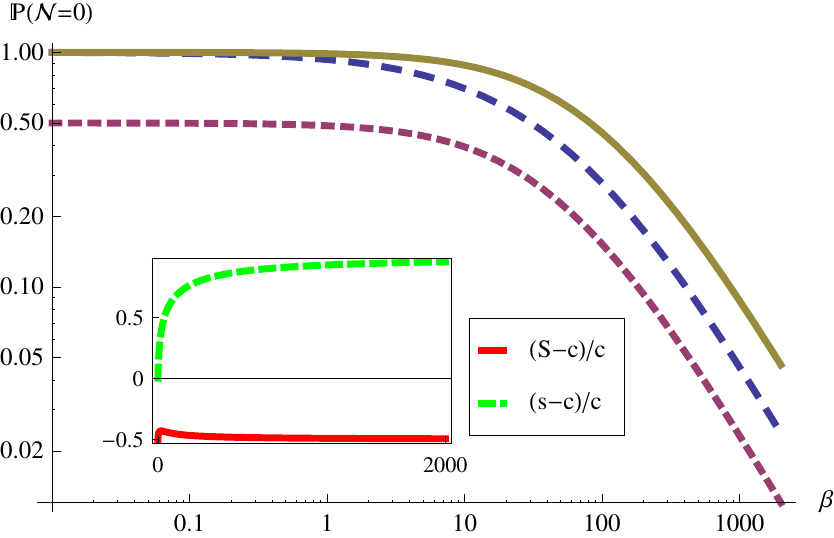}}
			\caption{\label{graph1}Probability of webpage ``death'' $\mathbb{P}(\mathcal{N}=0)$, in dependency on webpage rate constant $\beta$
			in critical ((a): $\lambda=\mu=5$; (b): $\lambda=\mu=50$), supercritical ((a): $\lambda=10$, $\mu=5$;
			(b): $\lambda=50$, $\mu=25$), and subcritical ((a): $\lambda=5$, $\mu=10$; (b): $\lambda=50$, $\mu=100$) regimes.
			In the insets, relative distance between supercritical and critical (solid line), and subcritical and critical
			regimes (dashed line) with respect to the critical regime.}
		\end{figure}

			Let us move now to the cases in which the linear birth-death process governing the appearance of in-links is supercritical or subcritical,
			i.e.,\ $\lambda>\mu$ or $\lambda<\mu$. For notational purposes let us define the functions (with $x,y,b,\nu>0$)
			\begin{align*}
				& r_b(x,y) = c_b(x,y)
				{}_2 F_1 \left( 1,\frac{b}{x-y}; 2+\frac{b}{x-y}; \frac{y}{x} \right), \\
				& q_b^\nu(x,y) = d_b^\nu(x,y) {}_2 F_1 \left( \nu+1, 1+\frac{b}{x-y};
				\nu+1+\frac{b}{x-y}; \frac{y}{x} \right),
			\end{align*}
			where ${}_2 F_1$ is the Gauss hypergeometric function and where
			\begin{align*}
				& c_b(x,y) = \frac{b}{x-y} \Gamma\left( \frac{b}{x-y} \right) \Bigr/ \Gamma
				\left( 2+\frac{b}{x-y}\right), \\
				& d_b^\nu(x,y) = \frac{b (x-y)}{x^2} \Gamma\left(1+\frac{b}{x-y}\right)\Gamma(\nu) \Bigr/
				\Gamma\left(\nu+1+\frac{b}{x-y}\right).
			\end{align*}
			By exploiting the integral representation
			\begin{align}
				\label{mattarello}
				{}_2 F_1 (a,b;c;z) = \frac{\Gamma(c)}{\Gamma(b) \Gamma(c-b)} \int_0^1 y^{b-1} \frac{(1-y)^{c-b-1}}{
				(1-yz)^{a}} \, \mathrm dy,
			\end{align}
			(see \cite{abramowitz}, page 558) considering \eqref{su}
			and after tedious calculations we obtain that in the supercritical case ($\lambda>\mu$)
			\begin{align}
				\mathbb{P}(\mathcal{N}=0) = (\mu / \lambda) r_\beta(\lambda,\mu),
			\end{align}
			and
			\begin{align}
				\label{dec}
				\mathbb{P}(\mathcal{N}=n) =
				q_\beta^n (\lambda,\mu),
			\end{align}
			while in the subcritical case ($\lambda<\mu$)
			\begin{align}
				\mathbb{P}(\mathcal{N}=0) = r_\beta(\mu,\lambda),
			\end{align}
			and
			\begin{align}
				\label{dec2}
				\mathbb{P}(\mathcal{N}=n) =
				(\lambda/\mu)^{n+1} q_\beta^n (\mu,\lambda).
			\end{align}
			The examples of the limiting distribution are depicted in Fig.\ \ref{figure2b}(a,c).
			Notably, in the supercritical case a power-law behavior is present, while in the subcritical case the tail of the distribution
			decays more rapidly. Moreover, if $\lambda < \mu$, the subcritical character of the process governing the arrival of in-links
			dominates the dynamics of the whole system and the parameter $\beta$ plays an almost irrelevant role.
			Fig. \ref{graph1} presents two examples of the $\mathbb{P}(\mathcal{N}=0)$,
			as a function of the webpage rate constant $\beta$, in the three cases, supercritical, critical, and subcritical.
			The insets highlight the behaviors
			of the relative distances between the three curves with respect to different choices of the parameters $\lambda$ and $\mu$.
			We retrieve the Yule--Simon distribution \eqref{yule} for fixed $\lambda>0$ and $\mu \to 0$ (classical Yule model)
			from \eqref{dec} by recalling
			that ${}_2 F_1 (a,b;c;0) = 1$. Clearly, in this case
			$\mathbb{P}(\mathcal{N}=0) = 0$.
			
			Manipulation of the probabilities is made simple by recalling that Gauss Hypergeometric functions admit a huge varieties of
			different representations and relations (see e.g.\ \cite{MR2360010}, Sections 9.1, 9.2,
			\cite{spanier}, Chapter 60, or \cite{MR2723248}, Chapters 13, 15). For example, when $\lambda>\mu$ and $n \ge 1$, by considering the representation
			15.2.1 of \cite{MR2723248} we have that $\mathbb{P}(\mathcal{N}=n)$ can be written also as the absolutely convergent series
			\begin{align}
				\mathbb{P}(\mathcal{N}=n) = \frac{\beta(\lambda-\mu)}{n\lambda^2} \sum_{r=0}^\infty
				\frac{\Gamma(n+1+r)\Gamma\left(1+r+\frac{\beta}{\lambda-\mu}\right)}{
				r! \Gamma\left(n+1+\frac{\beta}{\lambda-\mu}+r\right)} \left(\frac{\mu}{\lambda}\right)^r,
			\end{align}
			or by means of generalized Wright functions ${}_2\Psi_1$ as
			\begin{equation*}
				\mathbb{P}(\mathcal{N}=n) = \frac{\beta(\lambda-\mu)}{n\lambda^2} {}_2\Psi_1
				\left[ \left.
				\begin{array}{l}
					(n+1,1)\left(1+\frac{\beta}{\lambda-\mu},1\right) \\
					\left( n+1+\frac{\beta}{\lambda-\mu},1 \right)
				\end{array}
				\right| \frac{\mu}{\lambda}
				\right].
			\end{equation*}
			The representation of the probabilities by means of the power series expansion of the Gauss hypergeometric function provides a simple and efficient
			means for its computation (see \cite{temme}, Section 2.3.1. 
			Similar representations hold for all the other cases $\lambda=\mu$ and $\lambda < \mu$.

	\section{Generalized model: moments}		
		
		\label{mome}
		To complete the analysis of the model we calculate the moments of $\mathcal{N}$.
		In order to derive the expectation $\mathbb{E} \, \mathcal{N}$ it is
		useful to realize that it can be written as
		\begin{align*}
			\mathbb{E} \, \mathcal{N} = \int_0^\infty \beta e^{-\beta t} \mathbb{E} N_{\lambda,\mu}(t) \, \mathrm dt.
		\end{align*}
		Therefore we immediately have
		\begin{align}
			\label{media}
			\mathbb{E} \, \mathcal{N} & = \beta \int_0^\infty e^{-\beta t} e^{(\lambda-\mu)t} \mathrm dt
			= \frac{\beta}{\beta - (\lambda-\mu)}.
		\end{align}
		The above integral converges for $\mu>\lambda$ and, if $\lambda > \mu$, for $\beta > (\lambda-\mu)$,
		otherwise $\mathbb{E} \, \mathcal{N} = \infty$.
		It holds $\mathbb{E} \, \mathcal{N} = 1$ for $\lambda=\mu$.
		
		To derive the variance, we first determine the second moment.
		By recalling that for the classical linear birth-death process
		\begin{align*}
			\mathbb{E}\left[ N_{\lambda,\mu}(t) \right]^2 = \frac{2\lambda}{\lambda-\mu} e^{2(\lambda-\mu)t}
			- \frac{\lambda+\mu}{\lambda-\mu}
			e^{(\lambda-\mu)t}, \qquad \lambda \ne \mu,
		\end{align*}
		we calculate the second moment of $\mathcal{N}$ as
		\begin{align*}
			\mathbb{E}\,\mathcal{N}^2 = \frac{2\lambda}{\lambda-\mu} \frac{\beta}{\beta-2(\lambda-\mu)}
			-\frac{\lambda+\mu}{\lambda-\mu} \frac{\beta}{\beta-(\lambda-\mu)},
		\end{align*}
		and therefore
		\begin{align}
			\label{var}
			\mathbb{V}\text{ar} \, \mathcal{N} = \frac{\frac{2 \lambda\beta}{\beta-2(\lambda-\mu)}
			-\frac{(\lambda+\mu)\beta}{\beta-(\lambda-\mu)}}{\lambda-\mu} - \left( \frac{\beta}{\beta-(\lambda-\mu)} \right)^2.
		\end{align}
		The integral giving rise to the above expressions converges
		if $\beta/(\lambda-\mu)>2$, otherwise $\mathbb{V}\text{ar}\, \mathcal{N}=\infty$.
		For $\mu=\lambda$, as $\mathbb{E}\left[ N_{\lambda,\mu}(t) \right]^2 = 2\lambda t +1$,
		\begin{align*}
			\mathbb{E} \, \mathcal{N}^2 = 1+2 \lambda \beta \int_0^\infty e^{-\beta t} t \, \mathrm dt
			= 1+2 (\lambda/\beta). 
		\end{align*}
		It immediately follows that
		\begin{align}
			\mathbb{V}\text{ar} \, \mathcal{N} = 2(\lambda/\beta).
		\end{align}

		Higher-order moments can be derived from the probability generating function
		$\mathcal{G}(u) = \sum_{n=0}^\infty u^n \mathbb{P}(\mathcal{N}=n)$. This is determined
		by similar reasoning as before but this time a different special function is involved.
		In particular the so-called Appell hypergeometric function (see \cite{erdelyi}, Chapter 5) will be used,
		\begin{align}
			\label{appell}
			& F_1(a;b_1,b_2;c;z_1,z_2) \\ & = \frac{\Gamma(c)}{\Gamma(a)\Gamma(c-a)} \int_0^1 \frac{y^{a-1} (1-y)^{c-a-1}}{
			(1-y z_1)^{b_1} (1-yz_2)^{b_2}} \, \mathrm dy. \notag
		\end{align}
		Let $G(u,t)$ be the probability generating function of the linear birth-death process $N_{\lambda,\mu}(t)$.
		The probability generating function of $\mathcal{N}$ can be determined by using
		$\mathcal{G}(u)=\int_0^\infty \beta e^{-\beta t}G(u,t) \, \mathrm dt$,
		and \eqref{appell}. After some computations, for $\lambda > \mu > 0$, we get
		\begin{align}
			\mathcal{G}(u) = (\mu/\lambda) f_\beta \left( \lambda,\mu,\frac{\lambda u-\mu}{\mu(u-1)}, \frac{\lambda u-\mu}{
			\lambda(u-1)} \right), \quad |u|\le 1,
		\end{align}
		where
		\begin{align*}
			f_b(x,y,z,s) = F_1 \left( \frac{b}{x-y};-1,1;1+\frac{b}{x-y};z,s \right).
		\end{align*}
		For the subcritical case, i.e.\ when $\lambda<\mu$, we obtain
		\begin{align}
			\mathcal{G}(u) = f_\beta \left( \mu,\lambda, \frac{\mu(u-1)}{\lambda u-\mu}, \frac{\lambda(u-1)}{\lambda u-\mu} \right),
			\quad |u|\le 1.
		\end{align}
		The critical case $\lambda=\mu$ is simpler to cope with and results in a form which involves a difference of confluent hypergeometric
		functions. In particular we have
		\begin{align}
			\mathcal{G}(u) = \frac{\beta \, u \, U\left(1,1,\frac{\beta}{\lambda(u-1)}\right)}{\lambda(u-1)} 
			- U\left( 1,0,\frac{\beta}{\lambda(u-1)} \right).
		\end{align}

		\section{Comparison with real data}	
		
			\label{ddataa}
			\begin{figure}
				\centering
				\includegraphics[scale=0.85]{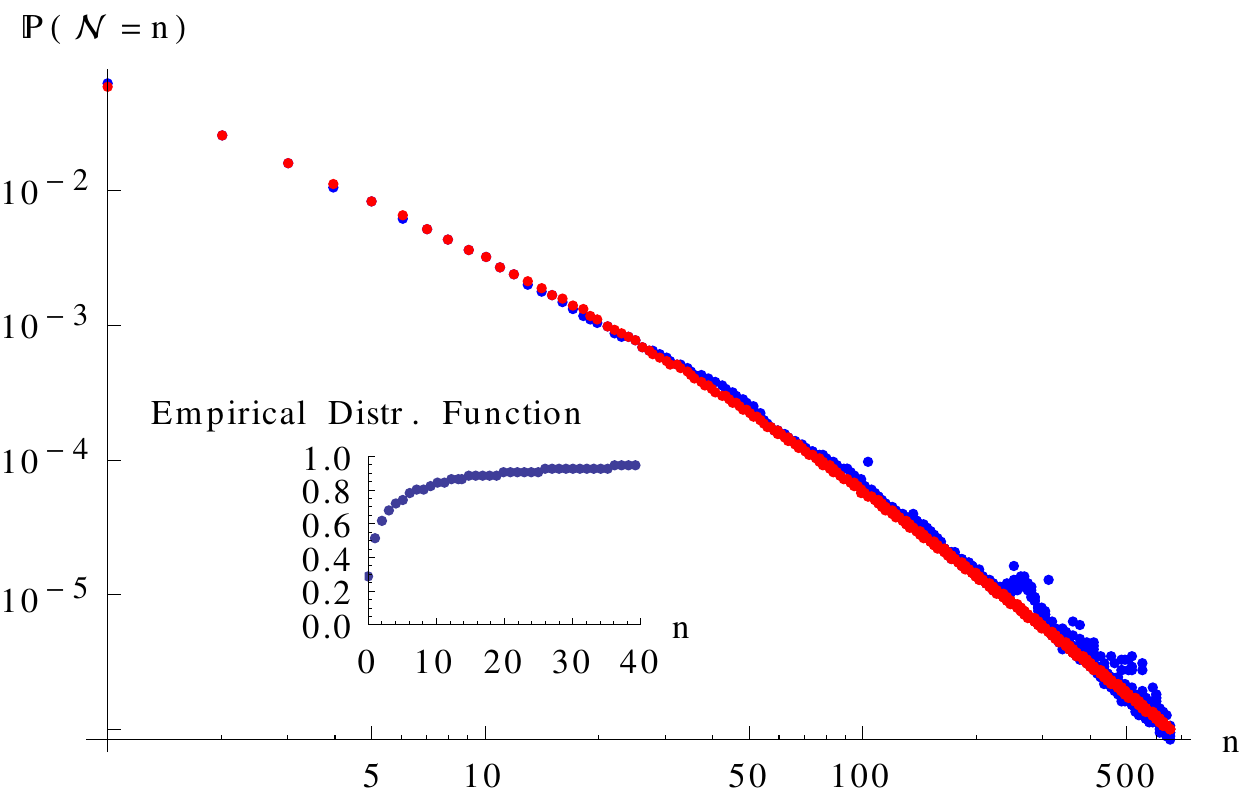}
				\caption{\label{data}The generalized Yule model compared to the empirical probability mass function
					for the number of in-links in the World Wide Web ($n \ge 1$). Data are taken from
					Web Data Commons, University of Manheim. It can be seen that the model (in red)
					fits the data even for small values of $n$ (with parameters $(\lambda;\mu;\beta) = (4;3.8;0.272)$).
					In the inset, the empirical distribution function calculated on the data. Note the importance of the small values of $n$ as
					the cumulative distribution function immediately saturates.}
			\end{figure}		
			As it is well known, data from World Wide Web show a scale free behavior. The proposed generalized model exhibits
			this feature in the supercritical regime.
			Having the complete distribution of in-links we are able to compare our model predictions against real data.
			The most recent available data on in-links of World Wide Web we found were collected in 2012 (Web Data Commons, University of
			Manheim\footnote{http://webdatacommons.org/hyperlinkgraph/}). In Fig.~\ref{data} we show the excellent performance
			of the proposed model with respect to these data. Note the goodness of fit also in the case of a few in-links.
			Furthermore, let us stress that most of the probability concentrates on points corresponding to a very few in-links
			(the cumulative distribution function at $n=20$ is already $0.9$) and this tells us that it is essential to derive closed form
			representations for these probabilities.

	\section{Generalized model: effect of link removal}
			
		\label{fuego}
\begin{figure}
			\centering
			\includegraphics[scale=1]{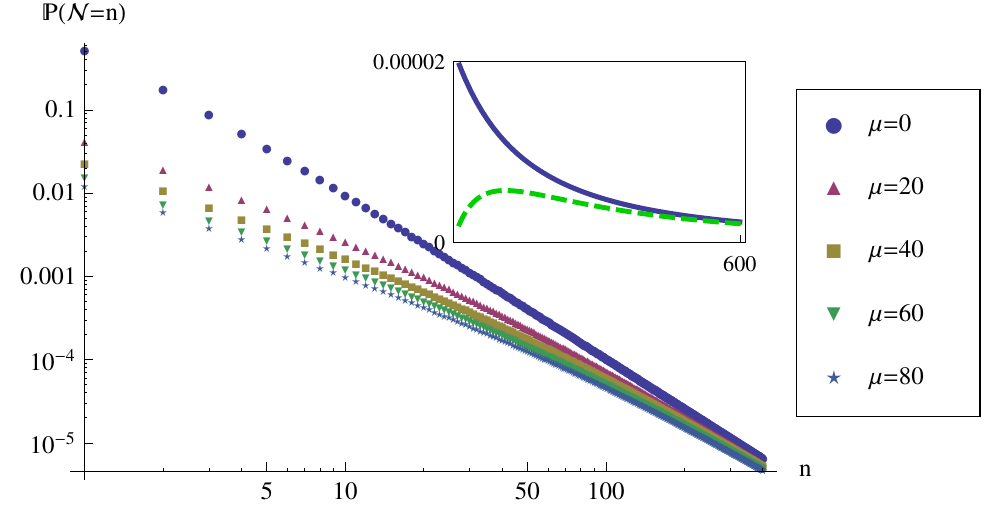}
			\caption{\label{new}The probability distribution $\mathbb{P}(\mathcal{N}=n)$, $n \ge 1$, for different values of
			removal rate $\mu$
			(note that $\mu = 0$ corresponds to Yule model) and
			$\beta=\delta=1$. Note that, with the choice $\beta=\delta$,
			the tail behavior for different values of $\mu$ can be compared more easily (see also Fig. \ref{cinque}
			for an example of $\beta\ne \delta$).
			In the inset (solid line) the theoretical distribution
			$\mathbb{P}(\mathcal{N}=n)$ for $\mu=80$, $\delta=1$, and (dashed line) its
			asymptotic behavior.}
		\end{figure}

\begin{figure}
			\centering
			\includegraphics[scale=1]{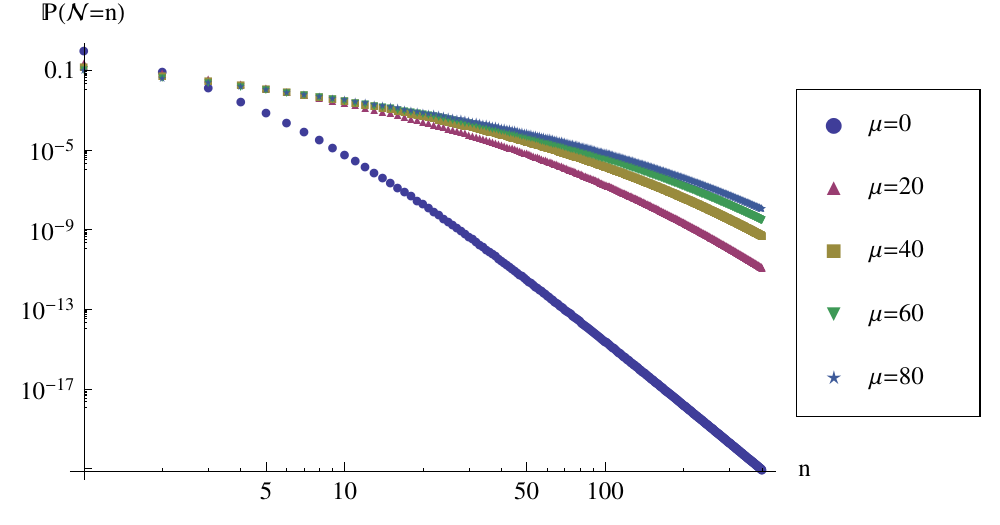}
			\caption{\label{cinque}The probability distribution $\mathbb{P}(\mathcal{N}=n)$, $n \ge 1$, for different values of $\mu$
			and $(\beta,\delta)=(10,1)$. Notice that, since $\delta=1$, the underlying behavior is supercritical but
			$\lambda$ increases together with $\mu$.}
		\end{figure}

		In the previous sections we presented a generalization of Yule model characterized by the introduction
		of the possibility of link detachment. This resulting model
		is more complex and one may wonder about the necessity of this added complexity.
		In the supercritical case, an intuition could suggest to replace the proposed model with a simple Yule model characterized by a decreased
		intensity $\delta = \lambda-\mu$ for the links generation. We show here that this intuition is misleading.
		The most evident difference is that only in the proposed model the event $\{\mathcal{N}=0\}$ has a positive probability.				
		A further evidence is highlighted by reparametrizing
		formula \eqref{dec} for $n \ge 1$ setting $\delta=\lambda-\mu$, $\lambda=\mu+\delta$, resulting in
		\begin{align*}
			\mathbb{P} (\mathcal{N}=n) = {} & \frac{\beta \delta \Gamma \left( 1+\frac{\beta}{\delta} \right) \Gamma(n)
			}{(\mu+\delta)^2\Gamma\left( n+1+\frac{\beta}{\delta} \right)} \\
			& \times {}_2 F_1 \left(n+1,1+\frac{\beta}{\delta};n+1+\frac{\beta}{\delta};\frac{\mu}{\mu+\delta}\right), \notag
		\end{align*}
		which can be written as
		\begin{align}
			\label{petr}
			\mathbb{P} (\mathcal{N}=n) = {} & \frac{\beta \Gamma \left( 1+\frac{\beta}{\delta} \right) \Gamma(n)
			}{\delta \Gamma\left( n+1+\frac{\beta}{\delta} \right)} \left( \frac{\delta}{\mu+\delta} \right)^{1-\frac{\beta}{\delta}} \\
			& \times {}_2 F_1 \left(\frac{\beta}{\delta},1+\frac{\beta}{\delta};n+1+\frac{\beta}{\delta};-\frac{\mu}{\delta}\right). \notag 
		\end{align}
		The last step comes from formula (1.3) of \cite{olde}.
		Note that the first factor of \eqref{petr} coincides with the analogous probability for the
		Yule model of parameter $\delta$. However, the second and third factors depend both on $\delta$ and $\mu$.
		The same feature characterizes the variance \eqref{var} while the mean \eqref{media} depends only on $\delta$. 
		Figs.~\ref{new} and \ref{cinque} show the shape of the distribution of $\mathcal{N}$ for different values of the death intensity
		$\mu$. We can see in Fig.\ \ref{new} how the initial part of the distribution, as $\mu$ becomes larger and larger, is
		more and more distant from a pure power-law behavior.

		An interesting property concerns the asymptotics of $\mathbb{P}(\mathcal{N}=n)$ for large values of $n$, i.e.\ the tail
		of the distribution. To investigate this aspect let us consider the asymptotics of the ratio
		\begin{align}
			\frac{\mathbb{P}(\mathcal{N}=n)}{\mathbb{P}(\mathcal{N}^Y=n)}
			= \left( \frac{\delta}{\mu+\delta} \right)^{1-\frac{\beta}{\delta}}
			{}_2 F_1 \left( \frac{\beta}{\delta},1+\frac{\beta}{\delta};n+1+\frac{\beta}{\delta};-\frac{\mu}{\delta} \right), \notag
		\end{align}
		where $\mathbb{P}(\mathcal{N}^Y=n)$ represents the probability of interest, \eqref{yule},
		for a classical Yule model of parameters $(\beta,\delta)$.
		When $n$ is large we have (see \cite{spanier}, formulae 60:9:3, page 605, and 48:6:1, page 474)
		\begin{align}
			\frac{\mathbb{P}(\mathcal{N}=n)}{\mathbb{P}(\mathcal{N}^Y=n)} \sim
			\left( \frac{\delta}{\mu+\delta} \right)^{1-\frac{\beta}{\delta}}
			- \left( \frac{\delta}{\mu+\delta} \right)^{1-\frac{\beta}{\delta}}
			\frac{\frac{\beta}{\delta}\left( 1+\frac{\beta}{\delta} \right)\left( 1+\frac{\mu}{\delta} \right)}{
			n+1+\frac{\beta}{\delta}} + \mathcal{O}\left( n^{-2} \right). \notag
		\end{align}
		The reader should realize that the first order asymptotics of the distribution of $\mathcal{N}$ and $\mathcal{N}^Y$
		are characterized by the same power-law exponent. However, it holds
		\begin{align*}
			\left( \frac{\delta}{\mu+\delta} \right)^{1-\frac{\beta}{\delta}} \ne 1
		\end{align*}
		if $\mu>0$ (i.e.\ if removal of in-links is introduced) and $\beta\ne\delta$. Furthermore
		the second order contribution to the asymptotics depends upon $\mu$.
		The inset of Fig.~\ref{new} compares the distribution of $\mathcal{N}$ with its asymptotics.
				
		The asymptotics of $\mathbb{P}(\mathcal{N}=n)$ can be exploited to estimate the value of $\beta/\delta$ and $\mu/\delta$.
		Indeed, using exact asymptotics of hypergeometric functions,
		we consider the dominant term of $\mathbb{P}(\mathcal{N}=n)$ as
		\begin{align*}
			\mathbb{P}(\mathcal{N}=n) \sim n^{-1-\frac{\beta}{\delta}} \Gamma\left( 1+\frac{\beta}{\delta}
			\right) \frac{\beta}{\delta} 
			\left( \frac{\delta}{\mu+\delta} \right)^{1-\frac{\beta}{\delta}} + \mathcal{O}\left( n^{-2-\beta/\delta} \right)
		\end{align*}
		(see also \cite{reed} for a different derivation).
		Taking the logarithm and setting $y=\log n$ we get
		\begin{align}
			\label{e130}
			\log \, & \mathbb{P}(\mathcal{N}=n) \sim -\left( 1+\frac{\beta}{\delta} \right) y \\
			& + \left( 1-\frac{\beta}{\delta} \right) \log \left( \frac{1}{1+\mu/\delta} \right)
			+ \log \Gamma\left( 1+\frac{\beta}{\delta}
			\right) + \mathcal{O}\left( n^{-1} \right). \notag
		\end{align}
		An estimate of the slope and intercept of the straight-line in \eqref{e130} allows the estimation
		of $\beta/\delta$ and $\mu/\delta$. Hence, from data a possible non-zero value of $\mu$ can be detected.
		Similarly, statistical test on the intercept and the constant term of the regression equation can be performed with
		classical regression methods, determining the p-value of the considered tests.
		These test can then be adapted to perform test on $\mu/\delta$ and $\beta/\delta$.
		We do not detail further this fact because other statistical features could be investigated and this will be the subject of a future work.
		
		Another possibility to obtain from data information on the value of the
		death intensity $\mu$ would be for example by means of the behavior of the distribution of $\mathcal{N}$ for small values
		of $n$.
		
	\section{Discussion}
	
		Traditional models disregard the detachment phenomenon but it may play a relevant role
		in network dynamics.
		As it is remarked in \cite{barab}, the observed growth of the World Wide Web
		measured as the number of connected webpages in time is exponential.
		This suggests that the assumption of the Yule process for the appearing of new webpages is reasonable.
		However, the original model assumes that the in-links exist forever. To avoid this hypothesis
		we generalized the model by using linear birth-death processes in substitution of
		classical birth processes for in-links dynamics.
		The proposed model admits the presence of \emph{lost} webpages, that is of webpages
		with zero in-links. The growth of the number of in-links per webpage is due to a linear birth
		phenomenon and therefore webpages with zero in-links are lost forever.		
		This feature of the model may be questionable and further modifications can be proposed.
	
		The introduction of removal was already considered in models motivated by different applications
		\cite{pnas,plos,reed,biol,patzkowsky,shi}	
		but here we obtain closed form results not yet appeared in the literature. In \cite{pnas,plos},
		the framework is that of macroevolutionary dynamics. Two different models are described and analyzed:
		the SEO model (speciation-extinction-origination model) and the BDM model (birth-death-mutation model) also studied in \cite{biol}.
		However, the microscopic dynamics of our model is different from all these models.
		In the SEO model, the origination rate for new families (new genera) is proportional to the number of species already present in the
		system while the BDM model is a discrete-time generalization of the Yule--Simon model allowing for deaths.
		The model discussed in \cite{patzkowsky} considers an origination rate for families which depends upon the number of species already
		present. Although an analytical treatment is performed, the overall behavior of the model is evaluated trough Monte-Carlo
		simulations.
		
		The framework of \cite{shi} is that of network growth. The model introduced and analyzed is
		a generalization of the Barab\'asi--Albert model in which removal of links is permitted. The model implements the
		usual preferential attachment for the link addition while an anti-preferential detachment is used for the link removal.
		The analysis is then performed by means of mean-field theory.		
		Other papers studied growth models allowing also negative fluctuations \cite{gabaix,saichev,pop}.
		These are works considering growth models with continuous state-space. They are considerably different from our
		work in the motivations, in the implemented models and also in the mathematical treatment.
	
		Finally, in the model treated in \cite{reed} (see also the references therein) the process governing the origination
		of new families is independent of that describing the growth of the number of species per family.
		Despite the different motivations and framework, the model considered in that paper coincides with that shown in the present paper.
		However, while \cite{reed} implements numerical and asymptotic methods, we arrive at closed form analytical results.

		We explicitely underline that differently from the Barab\'asi--Albert model, the Yule model disregards out-links.
		Their introduction would request further hypotheses on the model that are not necessary when the focus is on the popularity of a webpage,
		a quantity of interest for many applications.

		The results obtained in the present paper are rather general and therefore relevant to many different fields.
		Nonetheless, we wanted to stress the fact that network growth models allowing fluctuations in the number of in-links can
		be successfully modeled by means of suitable generalizations of the Yule model.

	\section{Concluding remarks}

		The main advantage of the investigated model lies in its mathematically tractability.
		We are able to present closed form expressions
		for the limiting distribution of the number of in-links per webpage. Its tail behavior is
		also analytically determined.
		The use of these formulae permits to recognize features not exhibited by the classical Yule model and to reproduce the observed
		behavior of in-links of World Wide Web.
		The main difference with the classical Yule model is that it did not allow the extinction of the in-links
		to a specific webpage. The presence of a non zero probability for $\{\mathcal{N}=0\}$ implies
		a decrease of the probabilities of $\{\mathcal{N}=n\}$, $n \ge 1$, and
		the probability of having a positive number of in-links is strictly smaller
		than the analogous quantity of the Yule model.
		This is also well illustrated by the fact that for the classical model always holds
		$\mathbb{E}\,\mathcal{N}^Y>1$ while this is not valid in the current version of the model.
		The comparison of the classical and generalized versions has been performed in the supercritical case, i.e.\ in the most relevant instance.
		The first-order tail behavior of the two models is similar but with a different proportionality
		coefficient. Differences arise in the second-order term.

		This is our first attempt to deal with the detachment issue and more realistic models can be considered.
		Growth phenomena admitting saturation or with more general nonlinear evolution rates are possible extensions.
		For example, a different situation arises if for a lost webpage there is a (possibly very small) probability of
		reobtaining an in-link. This circumstance would correspond to a birth-death process with immigration
		governing the increase of the number of in-links.
		Similarly, detachment can happen with higher probability when the number of in-links
		is low.			
		A model with removal of webpages in spite of existing in-links was considered in \cite{fenner} and could imply
		a variant of the studied model.

		We finally underline that our study is motivated by Web applications but scale-free networks are of interest in a
		variety of different contexts.
		Power-laws are the common feature of these models \cite{newman,barba,mit,hughes}.
		Social networks, macroevolutionary trees, scientific collaborations,
		words occurrences, neural networks, are some of examples of networks studied in the literature
		\cite{simon,pnas,plos,reed,patzkowsky,noh,grabowski,lansky}.

	\ack{The authors have been supported by project AMALFI (Universit\`{a} di Torino/Compagnia di San Paolo).
		Petr Lansky has been supported by Institute of Physiology RVO67985823 and the Centre for
		Neuroscience P304/12/G069.
	
		The authors would like to thank the anonymous referees for their useful comments and criticisms which improved the
		quality of the paper.}

\end{document}